\documentclass[fleqn,10pt]{wlscirep}
\usepackage[utf8]{inputenc}
\usepackage[T1]{fontenc}
\graphicspath{{images/}}
\usepackage{comment}
\usepackage{ragged2e}
\usepackage{algorithm,algorithmic}
\usepackage{siunitx}
\usepackage{rotating}

\title{RESIRE: real space iterative reconstruction engine for Tomography}

\author[1,*]{Minh Pham}
\author[2]{Yakun Yuan}
\author[2]{Arjun Rana}
\author[2]{Jianwei Miao}
\author[1]{Stanley Osher}
\affil[1]{Department of Mathematics, University of California, Los Angeles, CA 90095, USA}
\affil[2]{Department of Physics \& Astronomy and California NanoSystems Institute, University of California, Los
Angeles, CA 90095, USA}

\affil[*]{minhrose@ucla.edu}

\begin{abstract}
Tomography has made a revolutionary impact on diverse fields, ranging from macro-/mesoscopic scale studies in biology, radiology, plasma physics to the characterization of 3D atomic structure in material science. The fundamental of tomography is to reconstruct a 3D object from a set of 2D projections. To solve the tomography problem, many algorithms have been developed. Among them are methods using transformation technique such as computed tomography (CT) based on Radon transform and Generalized Fourier iterative reconstruction (GENFIRE) based on Fourier slice theorem (FST), and direct methods such as Simultaneous Iterative Reconstruction Technique (SIRT) and Simultaneous Algebraic Reconstruction Technique (SART) using gradient descent and algebra technique. In this paper, we propose a hybrid gradient descent to solve the tomography problem by combining Fourier slice theorem and calculus of variations. By using simulated and experimental data, we show that the state-of-art RESIRE can produce more superior results than previous methods; the reconstructed objects have higher quality and smaller relative errors. More importantly, RESIRE can deal with partially blocked projections rigorously  where only part of projection information are provided while other methods fail. We anticipate RESIRE will not only improve the reconstruction quality in all existing tomographic applications, but also expand tomography method to a broad class of functional thin films.  
We expect RESIRE to find a broad applications across diverse disciplines.
\end{abstract}
\begin{document}

\flushbottom
\maketitle
%
%
\thispagestyle{empty}


\section*{Introduction}
Tomography has widespread applications in physical, biological and medical sciences. Electron tomography lately has become a revolution in high-resolution 3D imaging of physical and biological samples.  In the physical sciences, atomic electron tomography (AET) has been developed to determine the 3D atomic structure of crystal defects such as grain boundaries, anti-phase boundaries, stacking faults, dislocations, chemical order/disorder and point defects, and to precisely localize the 3D coordinates of individual atoms in materials without assuming crystallinity \cite{miaoaaf2157, scott2012, chen2013, xu2015, haberfehlner2015}. 
Furthermore, using the advanced X-ray sources, coherent diffractive imaging (CDI) methods\cite{miao1999} can be integrated to produce 3D quantitative imaging of thick samples with resolutions in the tens of nanometers\cite{miao2006, nishino2009, jiang2010, dierolf2010, jiang2013, miao2015, holler2017}.

Along with the evolution of X-ray and electron Tomography, many different reconstruction algorithms are also developed to meet the demand for providing high-resolution 3D images from a finite number of projections. 
The algorithms fall into two main categories: one-step method filtered back-projection (FBP)\cite{frank2006, kak2002, herman2009} and iterative reconstruction (IR)\cite{dudgeon}.
FBP, well known for its efficient computation, works well when there are many enough projections and no missing data. Its fundamental is based on the relationship between Fourier slice theorem (FST) and Radon transform\cite{bracewell1995, lim1990} where the inverse responds to the reconstructions. Since working on  polar coordinates, inverse Radon transform requires the knowledge of all tilt angles in order to obtain a good reconstruction. In the cases of low radiation dose and geometric restrictions which produce inadequate data, FBP fails to reconstruct faithful objects and suffers artifacts. 

Real space iterative methods such as algebraic reconstruction technique (ART)\cite{gordon1970} can partially compensate for the artifacts. Modifications of ART such as simultaneous algebraic reconstruction technique (SART)\cite{andersen1984} and simultaneous iterative reconstruction technique (SIRT)\cite{gilbert1972} show significant improvement compared to FBP. These algorithms are developed based on a least square optimization problem to minimize the difference between measured and calculated projections. Gradient descent, an iterative method, is employed to refine reconstructions. Simplicity, fast running-time, and parallelization are the advantages of these methods.

Lately, Fourier-based iterative algorithms\cite{yang2017, miao2005, oconnor2006}, using information in both real and Fourier space as constraints in implementations, dramatically improve the performance. Equal slope tomography (EST)\cite{miao2005}, an example of such algorithm, has been successfully applied in AET to reconstruct the 3D atomic structure of crystal defects in materials, which shows successful recovery of diffraction pattern in the missing wedge direction\cite{scott2012, miaoaaf2157, chen2013, xu2015}. Additionally, with significantly lower radiation dose, EST can still produce reconstructions comparable to modern medical CT techniques\cite{fahimian2010, zhao2012, fahimian2013}. However, the drawback of EST is the strict requirement that the tilt angles must follow equal slope increments along a single tilt axis, which limits its broader applications. 

The EST development is followed by GENFIRE\cite{yang2017, pryor2017}, a generalized Fourier iterative reconstruction algorithm, which has been reported to produce high-resolution 3D imaging from a limited number of 2D projections. Using FST, GENFIRE transforms measurement constraint (measured projections) on real space into Fourier constraint on reciprocal space. The problem then becomes finding a 3D image that satisfies spatial constraints (i.e. real, positivity and boundary condition) and Fourier constraint. Alternating projection method is exploited to solve this classical two-constraint minimization problem. In addition, GENFIRE requires interpolation\cite{franke1982,shepard1968} to obtain Fourier values on Cartesian grid, and oversampling\cite{miao98,miao2000} to improve accuracy of this gridding process. These procedures are major drawbacks of GENFIRE since they require more memory (due to oversampling) and cause numerical error (caused by interpolation).

In this research, we motivate ourselves to develop a highly accurate Tomography gradient descent method by solving the least square error (LSE) problem. In the light of calculus of variations\cite{aubert2006}, we derive the gradient of the ``sum of squared errors" (SSE).
FST and interpolations are employed to compute the gradient which involves ``forward projection" and ``back projection" as defined by ART\cite{andersen1984}.  
The 3D object is then reconstructed iteratively with gradient descent method. 
These crucial points institute our proposed real space iterative reconstruction engine (RESIRE).

\section*{Method}
\subsection*{RESIRE algorithm}
\subsubsection*{Gradient descent}
Tomography is formulated as a least square optimization problem, minimizing the sum of squared errors (SSE)
\begin{equation} \label{min1}
    \min_O \; \varepsilon(O)  = \frac{1}{2} \sum_{\theta} \| \Pi_{\theta} (O) - b_{\theta} \|^2_F
\end{equation}
where $\Pi_{\theta}$ is the linear projection operator with respect to the tilt angle $\theta$ and $b_{\theta}$ is the corresponding measured projection. 
In a simple case where there is only one rotation axis (either x,y, or z), tilt angles can be denoted by a single variable. In the case of multiple rotation axes, object orientation can be expressed using Euler angles, which consist of a set of three angles with respect to a fixed coordinate system. More details regarding this topic will be discussed in later part. 
The error metric $\varepsilon$ is a net sum of every $\varepsilon_{\theta}$ caused by each measure projection $b_{\theta}$. We get a more explicit form for each error $\varepsilon_{\theta}$
\begin{equation} \label{min2}
    \varepsilon_{\theta}(O) = \frac{1}{2} \sum_{x,y} | \Pi_{\theta}(O) \{x,y\} - b_{\theta} \{x,y\} |^2
\end{equation}

We compute the gradient of the error metric $\varepsilon$ w.r.t. the object $O$ in the sense of calculus of variation method\cite{aubert2006}. The following derivation is a simplified (discrete) version of a continuous problem. Furthermore, in order for the derivation to be logical, we must assume the object is continuous. For a voxel at spatial location $\{u,v,w\}$, the gradient reads:
\begin{equation} \label{eqn:grad}
\begin{aligned}
    \frac{\partial \varepsilon_{\theta} }{\partial O \{ u,v,w \} } &= \sum_{x,y} \Big( \Pi_{\theta} (O) \{x,y\} - b_{\theta} \{x,y\} \Big)  \frac{ \partial }{ \partial O \{ u,v,w \} }  \sum_z O \Big\{ R_{\theta} 
    \begin{bmatrix} x\\ y\\z   \end{bmatrix} 
    \Big\} \\
    & = \Pi_{\theta} (O) \{x,y\} - b_{\theta} \{x,y\} \qquad  \text{where } 
    \begin{bmatrix} u \\ v \\ w \end{bmatrix} = R_{\theta} 
    \begin{bmatrix} x \\ y \\ z \end{bmatrix} \quad \text{ for some } z
\end{aligned}    
\end{equation}
The first line of Eqn. \ref{eqn:grad} is obtained by vanilla chain rule while the second line uses an assumption that voxels of object are independent from each other, i.e. $ \partial O \{x,y,z \} / \partial O \{u,v,w \} = 1$ if $\{x,y,z\} = \{u,v,w\}$ and 0 otherwise. Using the fact that the transpose of a rotation matrix is its inverse, we derive the following transformation
\begin{equation} \label{eqn:affine}
    \begin{bmatrix} x\\y \end{bmatrix} = 
    \begin{bmatrix} R_{1,1} & R_{2,1} \\ R_{1,2} & R_{2,2} \end{bmatrix}  \begin{bmatrix} u\\v \end{bmatrix}
    + \begin{bmatrix} R_{3,1} \\ R_{3,2} \end{bmatrix} w
\end{equation}
where $R_{i,j}$ is the $(i,j)^{th}$ element of the rotation matrix $R_{\theta}$. Note that all the voxels $O\{u,v,w\}$ are on Cartesian grid (integer coordinates) while related $\{x,y\}$ are not. To understand the formula, we decompose the procedure into two steps and borrow terminology in ART for explanation.

Step 1: computing ``forward projection." The projection $\Pi_{\theta}(O)$ is computed from the current iterative object $O$ via FST. To improve the accuracy, the 3D object is padded with zeros before the Fourier transform. The oversampling ratio of this padding procedure is typically in range of $[2,4]$ providing a good trade-off between the accuracy and efficiency. We now apply FST, i.e. the 2D Fourier transform of a projection w.r.t tilt angle $\theta$ represents a 2D plane slicing through the origin of the 3D Fourier transform of that object. Taking the inverse Fourier transform of this 2D slice and drop the size back to the original projection size, we obtain the desired projection. This process is like the ``forward projection'' step in ART. However, ART presents the ``forward projection" as a matrix (vector) multiplication while RESIRE employs FST and interpolation. As a result, accuracy of the ``forward projection" is improved with RESIRE. 

Step 2: computing ``back projection." The measured projections are subtracted from the forward projections to obtain the difference $\Pi_{\theta}(O) - b_{\theta}$. Now the task is to compute the gradient $\Pi_{\theta}^T (\Pi_{\theta}(O) - b_{\theta} )$, i.e. apply $\Pi^T$ on $\Pi_{\theta}(O) - b_{\theta}$. Eqn. \ref{eqn:affine} shows how to ``back project" a 2D image to a 3D object. Specifically, each $(u,v)$ slice of the gradient (with $w$ varied) is a linear transformation of the difference. Further notice that those $(u,v)$ slices where $w\neq 0$ are translations of the $(u,v)$ slice at the origin ($w=0$). The amount of translation is exactly the zero order term $[R_{3,1} \;\; R_{3,2}]^T w$ in Eqn. \ref{eqn:affine}.

On the other hand, SIRT presents the ``back projection" as a row-normalized transpose of the ``forward projection"
\begin{equation} \label{sirt}
    O^{k+1} = O^k - t \sum_i \frac{ a_i O^k - b_{\theta,i} }{a_i a_i^T } a_i^T = O^k - t \hat{A}^T \Big(  A O^k - b_{\theta} \Big)
\end{equation}
where $A$ is the projection operator defined by ART, $\{a_i\}_i$ are its row vectors, and $\hat{A}$ is its row-normalized matrix, i.e. the row vector $\hat{a}_i$ of $\hat{A}$ is normalized as $\hat{a}_i = a_i / (a_i^T a_i )$. $b_{\theta,i}$ is the $i^{th}$ element of $b_{\theta}$. SIRT in fact uses vector multiplication form (the first equality in Eqn. \ref{sirt}) rather than matrix multiplication form (the second equality). This strategy allows SIRT to enforce parallel computing on CPU or GPU\cite{palenstijn2011}. However, the downside of SIRT is its limitation to 2D case (single tilt axis). When extending to the multiple-tilt-axis case, computing and storing matrix $A$ is infeasible. Furthermore, $A$ is defined differently than $\Pi_{\theta}$ as ART does not require the linear transformation as RESIRE does. Instead, ART can choose arbitrary equidistant points for its line integral approximation\cite{andersen1984}. The restricted linear transformation in RESIRE is the major difference between these two algorithms and explains why RESIRE outperforms the other.

To finish this part, we recap that RESIRE derives the gradient of the SSE in Eqns. \ref{min1} \& \ref{min2} through calculus of variation to improve the accuracy of the forward and back projections. Furthermore, interpolation plays an important role in our algorithm as it is required in both steps. It first helps to apply FST by computing 2D slices through the origin of a 3D Fourier transform in ``forward projection" and then compute linear transformations in ``back projection".


\subsubsection*{Step size analysis} This is the tricky part of the gradient descent method. We hope to approximate the Lipchitz constant $L$ of the gradient such that the following inequality holds
\begin{equation}
    \|\nabla \varepsilon(O_1) - \nabla \varepsilon(O_2) \| \le L \| O_1 - O_2 \| \quad \forall  \; O_1, \, O_2
\end{equation}
Hence, we can choose the step size to be $1/L$ for convergence guarantee. Assuming that the sum of intensity is conserved under projection operator and using the fact that $l_2$ norm is invariant under rotation, we can reduce the step size analysis to the simple case where tilt angle $\theta= 0$. We further simplify the analysis to 2D case under assumption of one tilt axis. Without loss of generality, we can assume the 2D object has dimension $N_z \times N_x$ and the projection is applied along z direction (the first dimension). Hence, the forward projection matrix is just a simple row vector $\Pi = {\bf 1}^{1\times N_z}$ where all elements are ones. Similarly, the ``back projection" operator $\Pi^T = {\bf 1}^{N_z \times 1}$ is a column vector of ones, i.e. $\Pi^T/N_z$ uniformly distributes 1D image back to 2D image. Then the composed operator $\Pi^T \Pi = {\bf 1}^{N_z \times N_z}$ is all-one matrix. We can easily check its $l_2$ norm $\|\Pi^T \Pi\| = N_z$. We now derive the Lipchitz constant 
\begin{equation}
    \|\nabla \varepsilon_{\theta}(O_1) - \nabla \varepsilon_{\theta}(O_2) \| = \| \Pi^T \Pi ( O_1 - O_2) \| \le \| \Pi^T \Pi \| \; \| O_1 - O_2 \| = N_z \, \| O_1 - O_2 \| 
\end{equation}
Since there are $n$ projections that contribute to the gradient, the lipschitz constant increases by a factor of $n$, i.e. the accumulated Lipchitz constant becomes $L = nN_z$. Then the gradient descent step size can be presented as $t/(nN_z)$ where $t\in (0,1]$ is normalized step size. When object is sparse, larger t can work. Our experimental results show that the algorithm still converges with $t=2$. This value of t is also selected for all experiments in this paper.
We finalize our gradient descent method by showing the governing formula. At the $(k+1)^{th}$ iteration, the update reads
\begin{equation}
    O^{k+1}\{u,v,w\} = O^k \{u,v,w\} - \frac{t}{n N_z} \sum_{\theta}\Big( \Pi_{\theta} (O^k) \{x,y\} - b_{\theta}\{x,y\} \Big)
\end{equation}
As mentioned above, tilt angles are single values in the case of one rotation axis. However, multiple-rotation-axis Tomography appears more often in practice. In this case, Euler angles, a set of three angles, are used to describe the orientation of a 3D object with respect to a fixed coordinate system. Hence a rotation can be presented as a sequence of rotations about standard x,y and z axes. We prefer the Euler convention ZYX to demonstrate object rotations, i.e. each tilt angle is expressed as a triple set $(\phi,\theta,\psi)$ which corresponds to rotations about z, y, and x axes in this respected order. The composed rotation matrix $R_{\{\phi,\theta,\psi\}} = Z_{\phi} Y_{\theta} X_{\psi}$ is a product of 3 basic rotation matrices $Z_{\phi}, \, Y_{\theta}$ and $ X_{\psi}$ where $Y_{\theta}$  representing a rotation matrix about Y axis by an angle $\theta$. Similar definitions are followed by $Z_{\phi}$ and $X_{\psi}$.

\subsubsection*{Pseudocode}
The method is summarized in algorithm \ref{alg1} as a pseudocode. Bi-linear interpolation is used in both ``forward projection" and ``back projection" due to its simplicity, efficiency, and fast running time.

\begin{algorithm}
\caption{RESIRE algorithm}
\label{alg1}
\textbf{Input}: set of N projections $\{b_{\theta_i}\}_{i=1}^N$ and corresponding tilt angles $\{\theta_i\}_{i=1}^N$,  number of iterations $K$, step size t \\
\textbf{Initialize}: $O^0$.
\begin{algorithmic}
\FOR {$k = 1,\dots,K$}
\FOR {$i=1,\dots,N$} 
\STATE compute ``forward projection" $\Pi_{\theta_i}(O^k)$ using FST and interpolation.
\STATE compute residual $\quad R_i (O^k) = \Pi_{\theta_i} (O^k) - b_{\theta_i}$
\STATE compute ``back projection" (or gradient) $\nabla \varepsilon_{\theta_i}(O^k)$ from $R_i(O^k)$ by formula $\ref{eqn:grad}$ using interpolation.
\ENDFOR \\
$O^{k+1} = O^k - \frac{t}{nN} \sum_{i=1}^N \nabla \varepsilon_{\theta_i}(O^k)$ 
\ENDFOR
\end{algorithmic}
\textbf{Output}: $O^K$
\end{algorithm}

\section*{Results}
\subsection*{Reconstruction from a numerical simulation: biological vesicle}

\begin{figure}
    \hspace{0.8cm} model \hspace{1.8cm} RESIRE \hspace{1.7cm} GENFIRE \hspace{1.8cm} FBP \hspace{2.2cm} SIRT \\
    \centering
    \includegraphics[width=3cm]{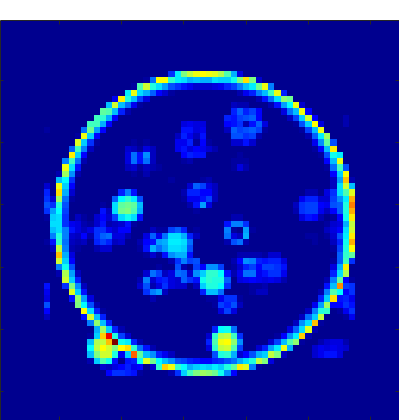}
    \includegraphics[width=3cm]{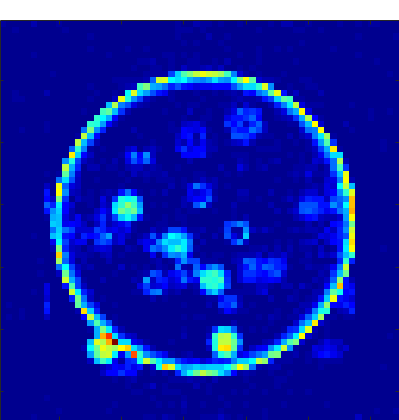}
    \includegraphics[width=3cm]{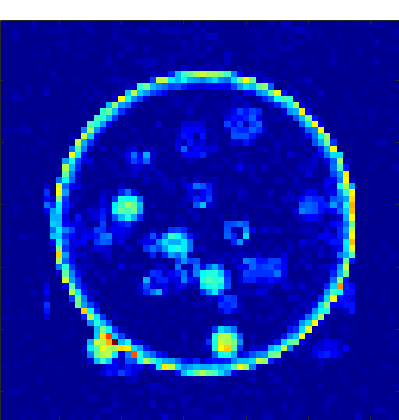}
    \includegraphics[width=3cm]{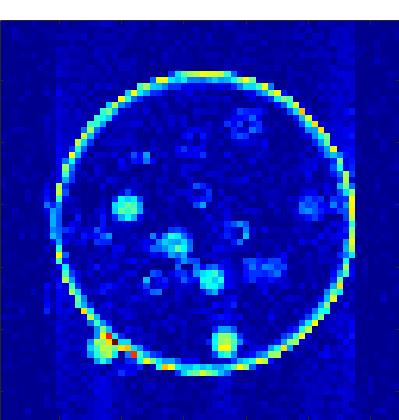}
    \includegraphics[width=3cm]{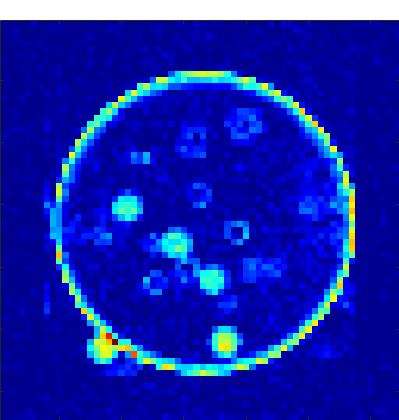}\\
    \includegraphics[width=3cm]{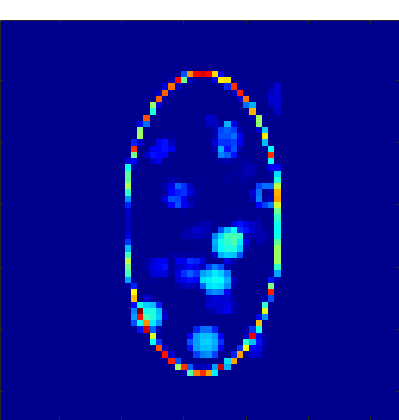}
    \includegraphics[width=3cm]{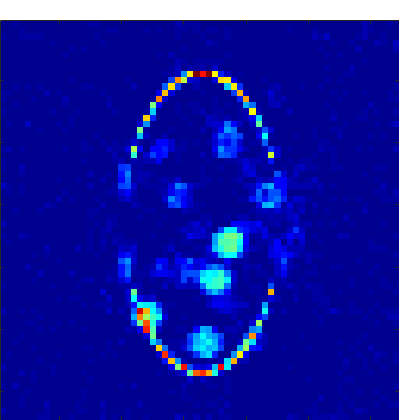}
    \includegraphics[width=3cm]{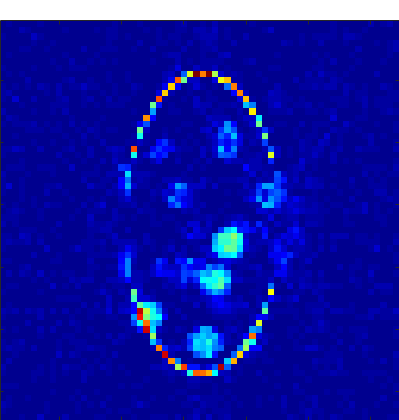}
    \includegraphics[width=3cm]{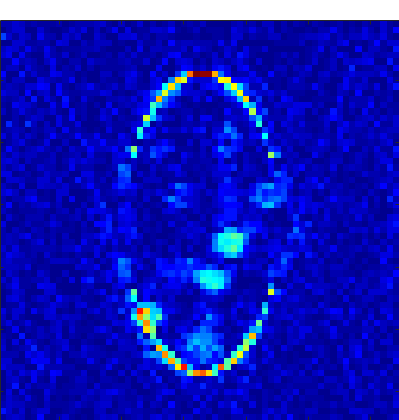}
    \includegraphics[width=3cm]{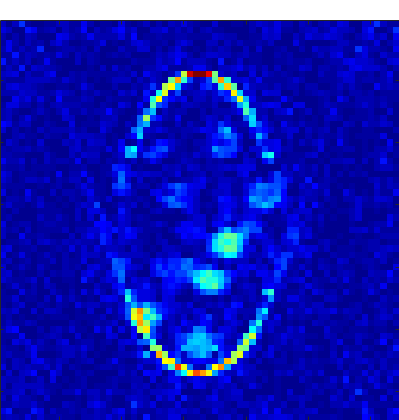} \\
    \includegraphics[width=3cm]{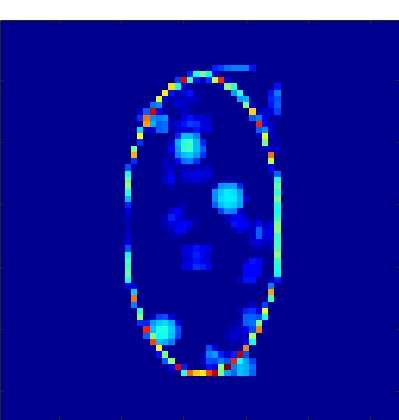}   
    \includegraphics[width=3cm]{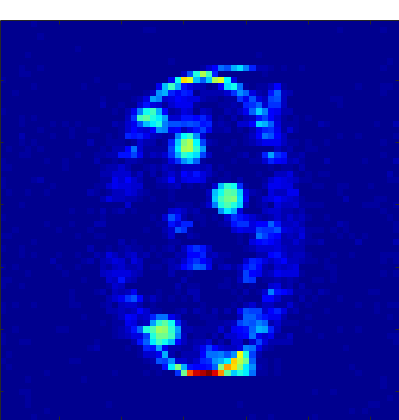} 
    \includegraphics[width=3cm]{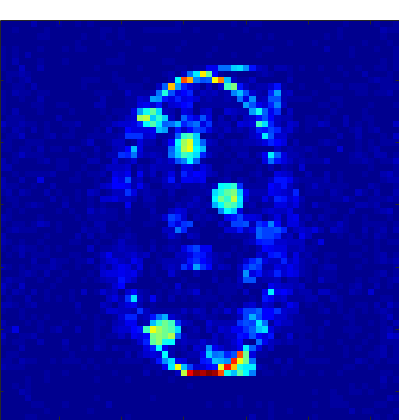} 
    \includegraphics[width=3cm]{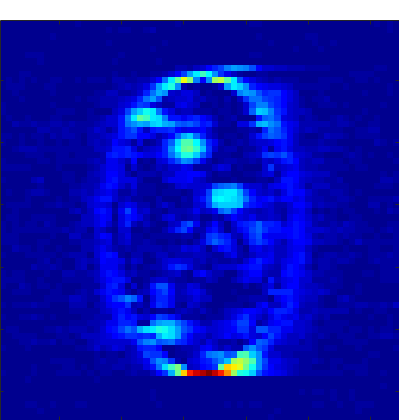} 
    \includegraphics[width=3cm]{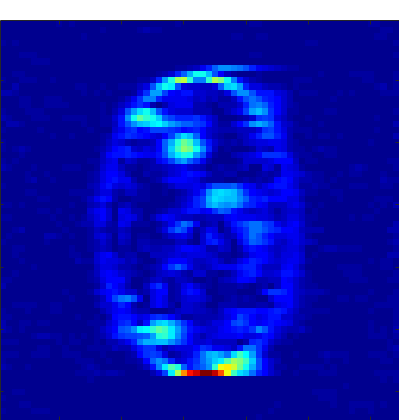} \\
    \includegraphics[width=8cm]{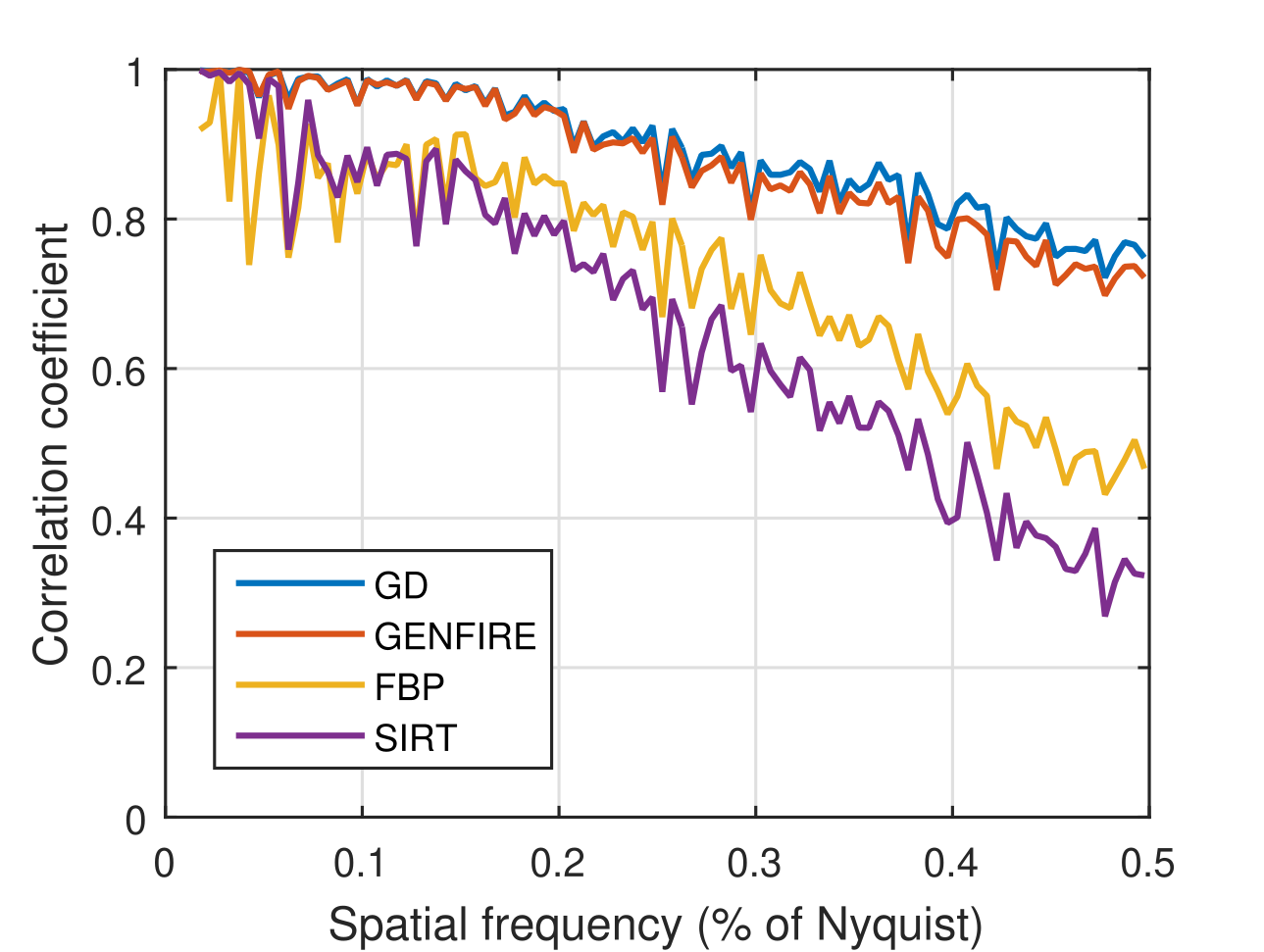}
    \caption{A vesicle model of $64 \times 64 \times 64$ voxel size from which 41 noisy projections simulated and 3D reconstructions using RESIRE, GENFIRE, FBP and SIRT. (a-c) Three 10-pixel-thick central slices of the vesicle model in the XY, XZ and YZ planes, respectively. The corresponding three reconstructed slices with RESIRE (d-f), GENFIRE (g-i), FBP (j-l), and SIRT (m-o), where the missing wedge axis is along the z-axis. The FSC between the reconstructions and the model shows that RESIRE produces a more faithful reconstruction than the other algorithms at all spatial frequencies.}
    \label{fig:vesicle}
\end{figure}

A model of biological vesicle of size $64\times 64 \times 64$ voxels was used in this test to demonstrate the performance of RESIRE algorithm. We assume $Y$ is the only rotation axis, then the other two Euler angles $\phi$ and $\psi$ are identically zeros.
Using FST, we calculate 41 projections whose tilt angles $\theta$ range from $-70^0$ to $70^0$ with step size $3.5^0$. Gaussian noise was added to projections to make to simulation data more realistic.
Reconstructions are performed using RESIRE, GENFIRE, FBP and SIRT. To monitor convergence, we use R-factor $R_F$ (relative error) as a metric to compare the relative difference between forward and measured projections.
\begin{equation}
    R_F = \frac{1}{n} \sum_{\theta} \; \Big[ \sum_{x,y}  |\Pi_{\theta} (O )\{x,y\} | - b_{\theta} \{x,y\} | \Big/ \sum_{x,y} | b_{\theta}\{x,y\} | \Big]
\end{equation}
All iterative reconstruction algorithms RESIRE, GENFIRE and SIRT were performed with 400 iterations. ASTRA Toolbox\cite{ vanAarle, van2015} was employed to achieve SIRT reconstruction while FBP reconstruction was performed using IMOD\cite{mastronarde1997}. The obtained R-factor are 28.59, 38.98, 35.44, 30.46$\%$ for RESIRE, GENFIRE, FBP and SIRT, respectively. Since directly minimizing the least square, RESIRE and SIRT obtain lower $R_F$ than the other methods. FBP obtains the largest $R_F$ because it is one step method. The $R_F$ of GENFIRE is also high due to the numerical error in the Fourier space interpolation.

Figure \ref{fig:vesicle} show 10-pixel-thick central slices of the 3D reconstructions in the XY (first row), XZ (second row) and YZ (third row) planes by RESIRE, GENFIRE, FBP and SIRT respectively, where the z-axis is the missing wedge direction.
Because there is no missing data in this direction, the XY central slices from all methods exhibit good agreement with the model.
However, along the missing wedge direction z, FBP and SIRT reconstructions show degradation, noise, and missing wedge artifact (shadow and smearing) while features obtained by RESIRE and GENFIRE appear to be finer and more isotropic. GENFIRE reduces but still shows some artificial noise in the XZ and YZ central slices because of the high numerical error in the Fourier space interpolation. On the other hand, RESIRE reduces significant noise, and obtains the finest and clearest images.


To quantify the results, we use Fourier shell correlation (FSC) to measures the normalized cross-correlation coefficient between the model and each reconstruction over corresponding shells in Fourier space (as a function of spatial frequency correlation). 
The FSC curves in Figure \ref{fig:vesicle} confirms that the RESIRE reconstruction is superior at all spatial frequencies compared to the other algorithms.

\subsection*{Reconstruction of a frozen hydrated cell}

\begin{figure}[ht]
    \hspace{1.8cm} RESIRE \hspace{3.cm} SIRT \hspace{3.2cm} FBP \smallskip \\
    \centering
    \includegraphics[width=4.cm]{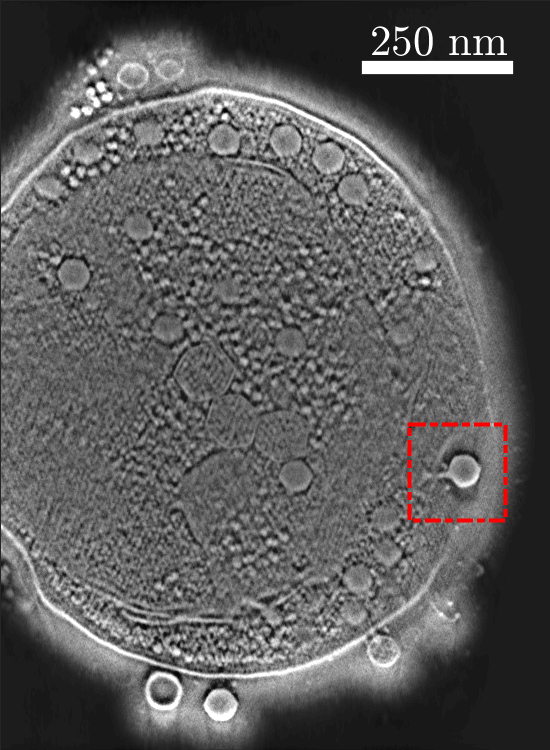}
    \includegraphics[width=4.cm]{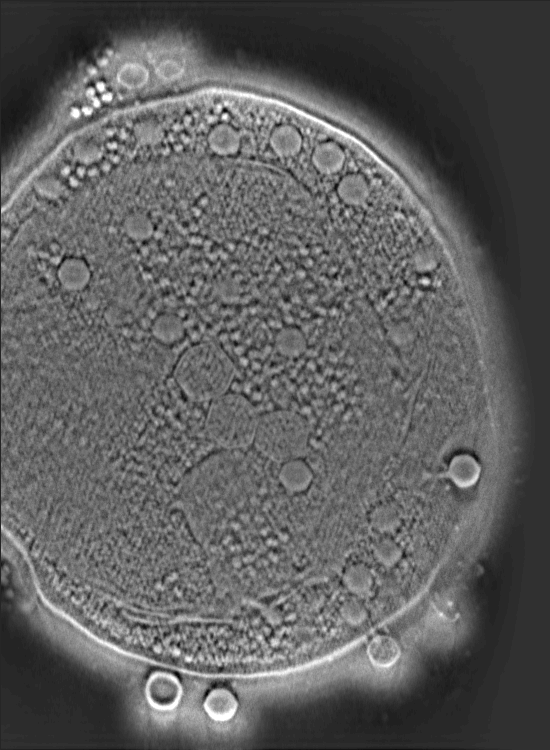}
    \includegraphics[width=4.cm]{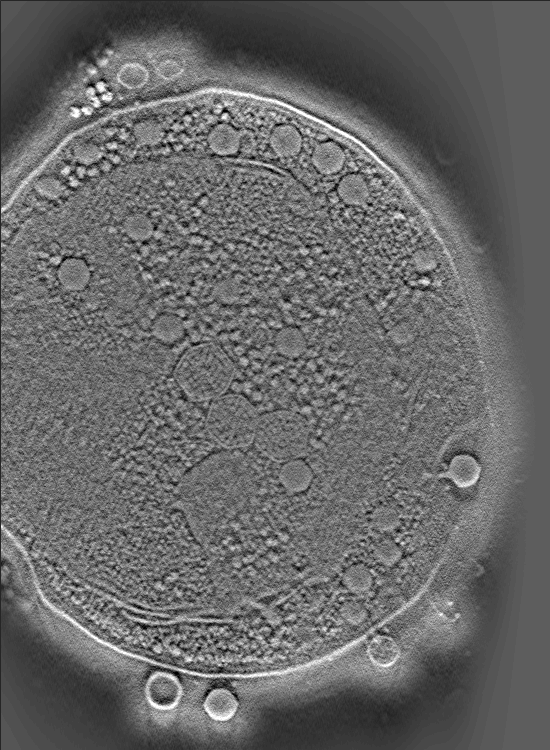} \\
    \includegraphics[width=4.cm]{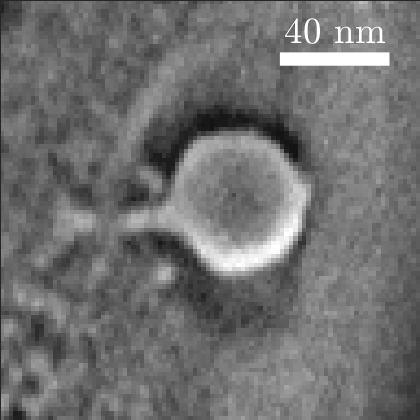}
    \includegraphics[width=4.cm]{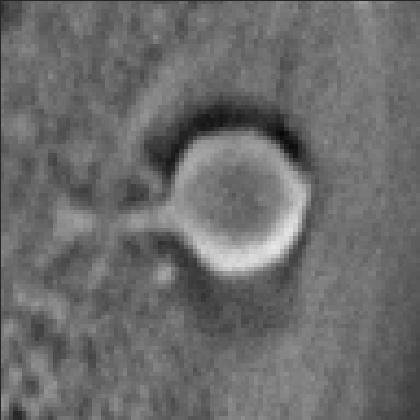}
    \includegraphics[width=4.cm]{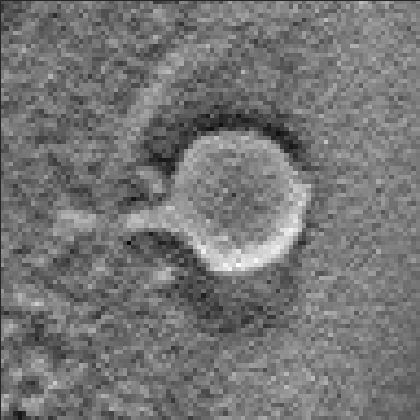} \\
    \caption{Reconstruction of a frozen-hydrated marine cyanobacterium, capturing the penetration of a cyanophage into the cell membrane. (a-c) 5.4-nm-thick (3-pixel-thick) slices of the cell in the XY plane reconstructed by RESIRE, SIRT, and FBP  respectively, where the y axis is the tilt axis and the beam direction is along the z axis. (d-e) Magnified views of the penetration of a cyanophage for the RESIRE, SIRT and FBP reconstructions in the XY, respectively. The high-quality projections results in good reconstructions. However, since RESIRE and SIRT use iterations, their reconstructions and R-factors are much better than FBP.  $R_F = 2.50\%, \, 4.39\%$ and $40.63\%$ for RESIRE, SIRT and FBP respectively. RESIRE obtains not only smaller $R_F$ but also better quality and resolution than SIRT and FBP. The membrane surrounding the cyanophage is more visible with RESIRE.}
    \label{fig:cyanophage}
\end{figure}

We run the reconstruction of a frozen-hydrated marine cyanobacterium in a late stage of infection by cyanophages. 42 projections were collected from a single tilt series, ranging from $-58^o$ to $+65^o$ with an increment of $3^o$. The projections were binned by $4 \times 4$ pixels and resulted in images with approximately $1.8 \times 1.8$ nm$^2$ per pixel. The details of experiment and pre-processing were explained in the GENFIRE paper\cite{pryor2017}. In this single tilt axis case, we can perform a comparison among RESIRE, SIRT and FBP algorithms. Both RESIRE and SIRT reconstructions were performed with 150 iterations.

Figure \ref{fig:cyanophage} shows a 5.4-nm-thick (3-pixel-thick) slice of the 3D reconstruction in the XY plane, capturing the penetration of a cyanophage into the cell membrane during the infection process. 
The SIRT and RESIRE reconstructions look similar and much better than FBP. 
The FBP reconstruction suffers peripheral noise through which features could not be seen clearly. This is because FBP is one-step method; hence, it cannot resolve the feature well with missing data.
On the other hand, both RESIRE and SIRT iterate to fit the measurements (measured projections). As a result, features in RESIRE and SIRT are resolved with higher quality.
These results are also consistent with R-factors where RESIRE (2.5\%) and SIRT (4.39\%) obtain much lower values than FBP (40.46\%). 
Because the projections have high quality and the noise level is low, both RESIRE and SIRT achieve good results. Therefore, their difference is marginal, and it is difficult to justify with the naked eye. Under careful comparison, the reconstruction by RESIRE can be verified to be better with sharper and clearer features, especially at the edges. In addition to this reasoning, the lowest R-factor reveals that RESIRE reconstruction fit the measurement better than other methods.

\subsection*{Reconstruction from an experimental data: Seriatopora aculeata coral skeleton}
\begin{figure}
    \hspace{2.cm} RESIRE \hspace{1.5cm} GENFIRE \hspace{2.cm} FBP \hspace{2.5cm} SIRT \\ 
    \centering
    \begin{turn}{90} \normalsize{XY projection} \end{turn}
    \raisebox{-.15\height}{ \includegraphics[height=3.4cm]{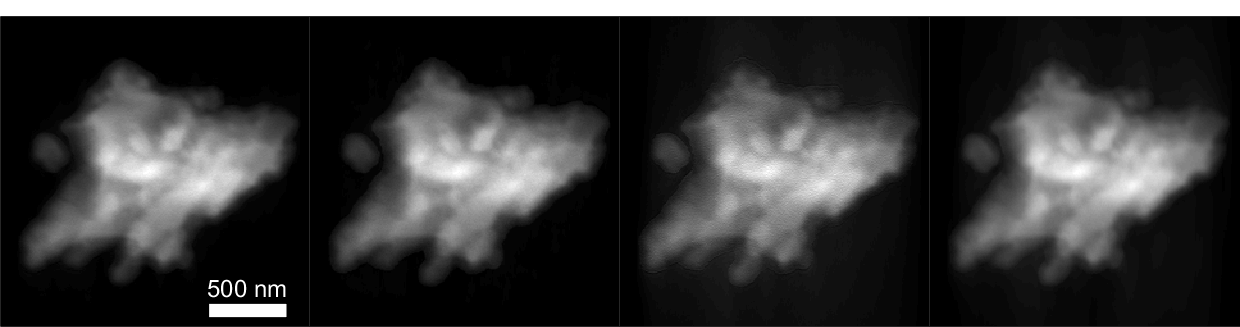} }\\
    \begin{turn}{90} \normalsize{XZ projection} \end{turn}
    \raisebox{-.15\height}{ \includegraphics[height=3.4cm]{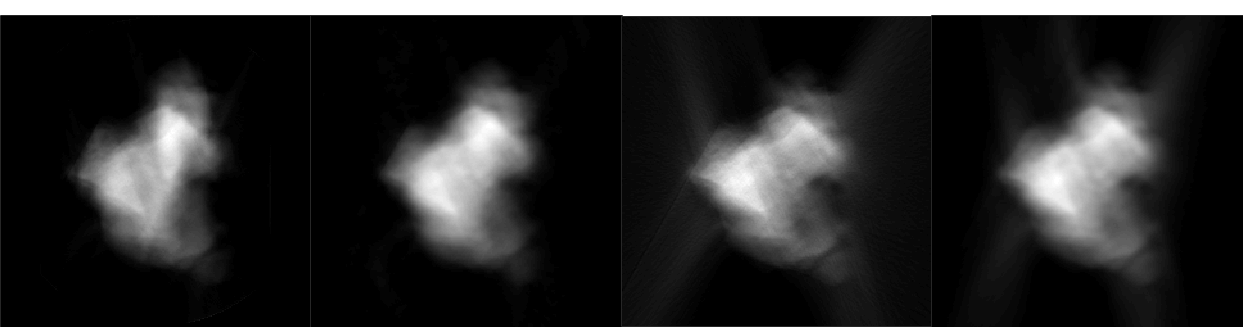}   }\\
    \begin{turn}{90} \normalsize{YZ projection} \end{turn}
    \raisebox{-.15\height}{ \includegraphics[height=3.4cm]{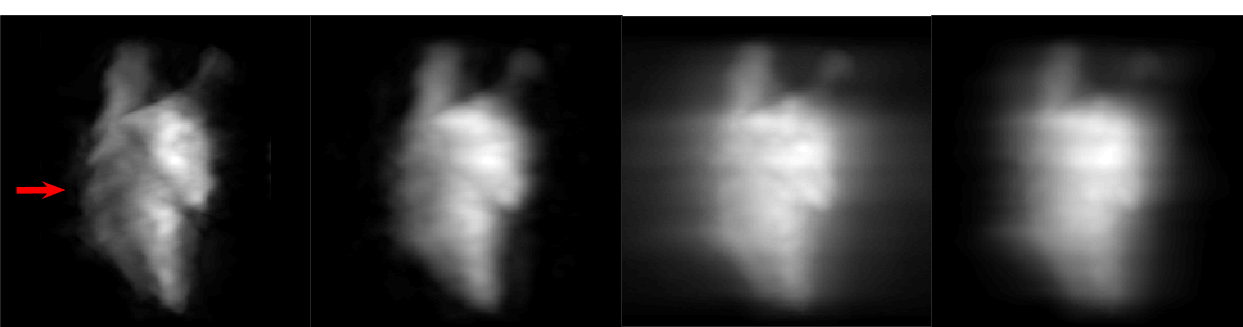} }\\
    \begin{turn}{90} \normalsize{XY central slices} \end{turn}
    \raisebox{-.15\height}{ \includegraphics[height=3.4cm]{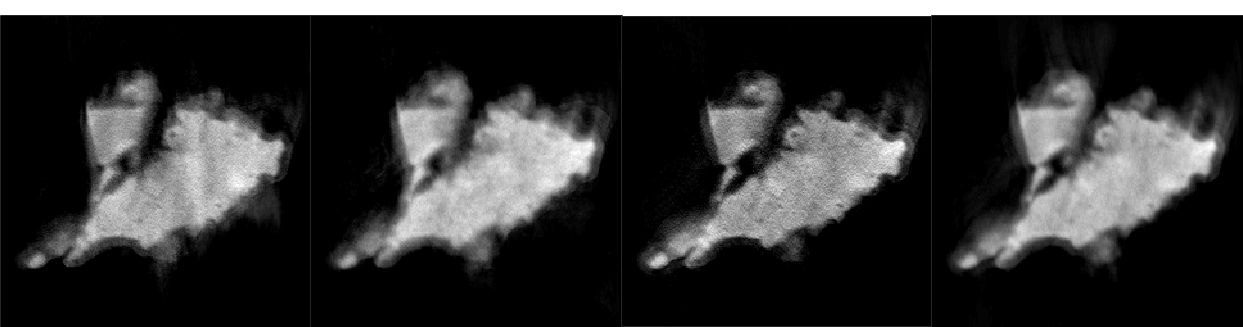} }\\
    \begin{turn}{90} \normalsize{XZ central slices} \end{turn}
    \raisebox{-.15\height}{ \includegraphics[height=3.4cm]{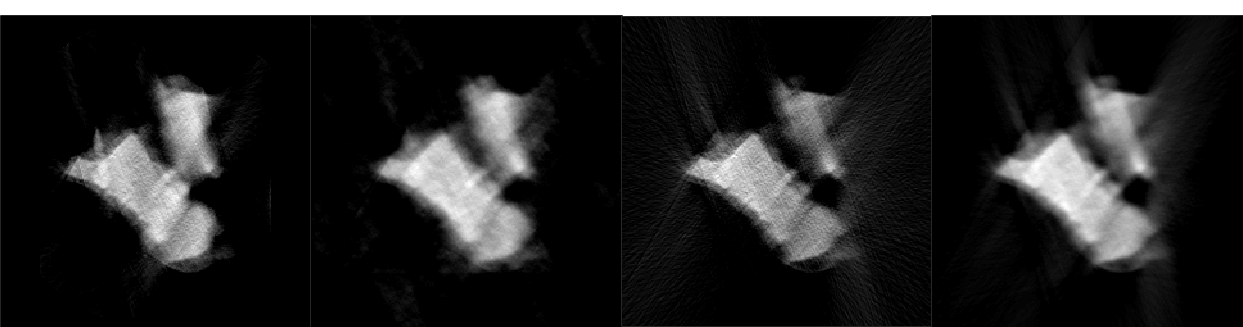} }\\
    \begin{turn}{90} \normalsize{YZ central slices} \end{turn}
    \raisebox{-.15\height}{ \includegraphics[height=3.4cm]{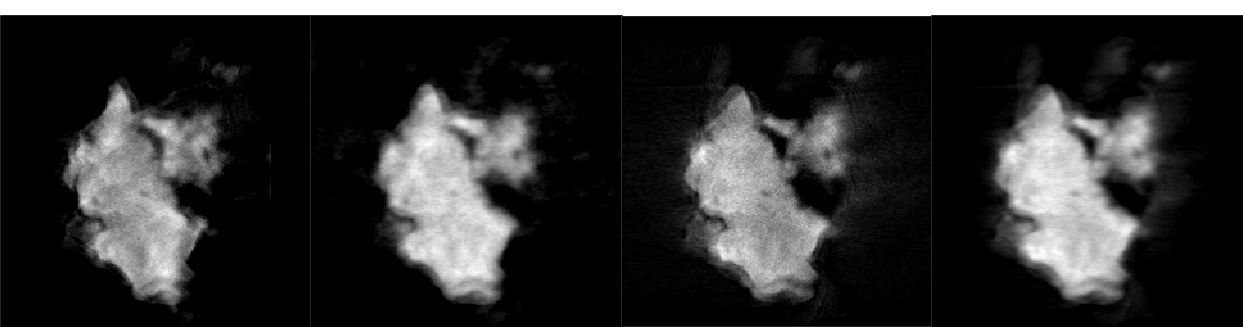} }
    \caption{Reconstruction of a Seriatopora aculeata coral skeleton by RESIRE, GENFIRE, FBP, and SIRT: projections and 30.3-nm-thick (3-pixel-thick) central slices in the XY, XZ and YZ planes, respectively. RESIRE reconstructs more faithful and finer structures inside the coral in all directions than other methods. Especially in the missing wedge direction z, while FBP and SIRT suffers blurring, elongation, peripheral noise, and line artifacts, RESIRE conquers these problems, and obtains the finest resolution. Furthermore, in the YZ projection features, indicated by a red arrow, can be only found visible with RESIRE.}
    \label{fig:coral}
\end{figure}

The Seriatopora aculeata sample used in this test is a pencil-thick, short, and tapered branch termed nubbin. Experiment was performed to obtain high-angle annular dark-field (HAADF-STEM) tomography data 
which includes a total of 69 images collected with a tilt range from -60 to 76 degrees in 2 degrees of increment. 
The projections were pre-processed with alignment, constant background subtraction and normalization.
Details of this experiment and pre-processing can be found elsewhere\cite{lo2019}.


After the pre-processing, Tomography reconstructions of the Seriatopora sample were performed by RESIRE, GENFIRE, FBP, and SIRT algorithms. Figure \ref{fig:coral} shows the reconstructed projections and 30.3-nm-thick (3-pixel-thick) central slices in XY, XZ and YZ planes. RESIRE is shown to completely outperform the other algorithms as the features obtained by RESIRE are finer, clearer, and more visible in all directions. Especially in the missing wedge direction, the YZ projections reconstructed by FBP and SIRT are corrupted by blurring and smearing. GENFIRE achieves a better reconstruction where significant artifact is reduced; however, the features are still blurred. RESIRE, in contrast, obtains the finest solution among these four algorithms. The features indicated by the red arrow in the YZ projection can be only obtained clearly and visibly with RESIRE.

\section*{Discussion}
In this article, we present RESIRE, a Tomography iterative reconstruction algorithm for high-resolution 3D objects. The method is designed to work with multiple rotation axes and is generally applicable to all Tomography applications. Through numerical simulations and experimental data on material science and biological specimens, RESIRE has been shown to produce superior reconstructions to several other most popular algorithms. More importantly, RESIRE can work well with limited projections and missing wedge while other methods fail and suffer artifacts such as elongation, blurring, and shearing. 

To summarize, RESIRE uses gradient descent to solve the least square problem $\|\Pi O - b\|^2$. Instead of pre-computing and storing the (forward) projection operator $\Pi$ as a matrix (such as SIRT does), RESIRE uses bi-linear interpolation and FST. 
Especially, FST helps to improve the forward projection while bi-linear interpolation saves enormous memory in the multiple-tilt-axis case where computing and storing operator matrix is extremely expensive. 
These two techniques are the major points that constitute our ``real space iterative reconstruction engine" RESIRE.

From mathematical point of view, ART, SIRT and SART approximate the forward projection (line integral) as a discrete sum. However, they do not specify how to choose grid points for this finite sum. RESIRE, in contrast, does linear transformation as detailed by Eqn. \ref{eqn:affine}, which significantly improves the performance  compared to SIRT. Moreover, RESIRE can exploit the refining technique, such as dividing a pixel into 4 sub-pixels used by Radon transform\cite{bracewell1995, lim1990}, to further improve the ``forward operator."

For global-transform-based methods such as GENFIRE, it converts measurements (projections) on real space to constraint on the reciprocal space via FST. This technique causes numerical error due to the gridding process (interpolation from non-grid points). Furthermore, by favoring smoothed solutions, reciprocal constraint prevents GENFIRE from obtaining high resolution object reconstruction. 
On the other hand, using a hybrid model (compute forward projection via FST and back projection via bi-linear interpolation), RESIRE can solve the downside issues of these methods. 
This explains why RESIRE has ability to work with high noise and produce more faithful reconstruction with higher resolution as well as finer and sharper features.

Furthermore, the extension of RESIRE to the thin film case allows us to work with non-isolated or extending samples which are highly popular in material science and industrial applications. Methods using global transform, such as FBP and GENFIRE, fails in this case because they interpret objects that are not consistently within the experimental FOV as noise. While, RESIRE can simply rigorously deal with this case.

Tomography is an ill-posed inverse problem and might need regularizer to compensate for the under-determination\cite{clason2019}. The optimal solution can be achieved after some iterations based on a suitable discrepancy principle\cite{lellmann2014}. This principle tells us a criterion where to stop the iterations according to prior noise information. The refinement will take place until the R-factor hits that certain value. Further iteration will improve the R-factor but not the result. Instead, the solution will move away from the optimal one.

Thanks to the simplicity, efficiency, and flexibility, RESIRE can have intensive applications in 3D image reconstructions, such as cryo-ET and vector tomography\cite{donnelly2017, donnelly2018, donnelly2018b} where reconstructing the 3D magnetic field of magnetic materials is a big interest. Looking forward, we expect RESIRE can be applied to a plethora of imaging modalities across different disciplines. 

\bibliography{sample}

\section*{Acknowledgements}
We thank Jihan for the preparation of the NiPt nanoparticle and Seriatopora aculeata samples. The HAADF-STEM imaging and ADF-STEM imaging with TEAM I was performed at the Molecular Foundry, Lawrence Berkeley National Laboratory, which is supported by the Office of Science, Office of Basic Energy Sciences of the US DOE under contract no. DE-AC02-05CH11231.
This work was also supported by STROBE: A National Science Foundation Science \& Technology Center under Grant No. DMR 1548924.

\section*{Author contributions statement}

J.M. directed the research; M.P. developed the main algorithm; M.P., Y.Y, A.R. develop the thin film extension; M.P performed the numerical simulation and the reconstruction of the experimental data. M.P., Y.Y. analyzed the results; M.P. wrote the manuscript with contributions from all authors. All authors reviewed the manuscript.

\section*{Additional information}
\textbf{Accession codes:} Matlab codes of the RESIRE algorithms and simulation data are available on our website.

\noindent \textbf{Competing Interests:} The authors declare that they have no competing interests.

\end{document}